\documentclass[12pt]{article}
\usepackage{amsfonts}
\usepackage{epsfig}
\usepackage{pdfpages}
\title{A Hyperbolic View of the Seven Circles Theorem}
\author{Kostiantyn Drach\thanks{\hskip 5pt Supported by the advanced grant ``HOLOGRAM'' of the European Research Council (ERC)}  \hskip 5pt and  Richard Evan Schwartz\thanks{\hskip 5 pt Supported by 
N.S.F. Research Grant DMS-1807320}}
\date{}

\newtheorem{theorem}{Theorem}[section]

\def\startproof{{\bf {\medskip}{\noindent}Proof: }}

\def\endproof{$\spadesuit$  \newline}

\def\C{\mbox{\boldmath{$C$}}}%
\def\H{\mbox{\boldmath{$H$}}}%

\begin{document}
\maketitle

\section{Introduction}

Cecil John Alvin Evelyn, or simply ``Jack'' to friends,
was born in 1904 in the United Kingdom in the aristocratic Evelyn family.

\begin{center}
\resizebox{!}{2in}{\includegraphics{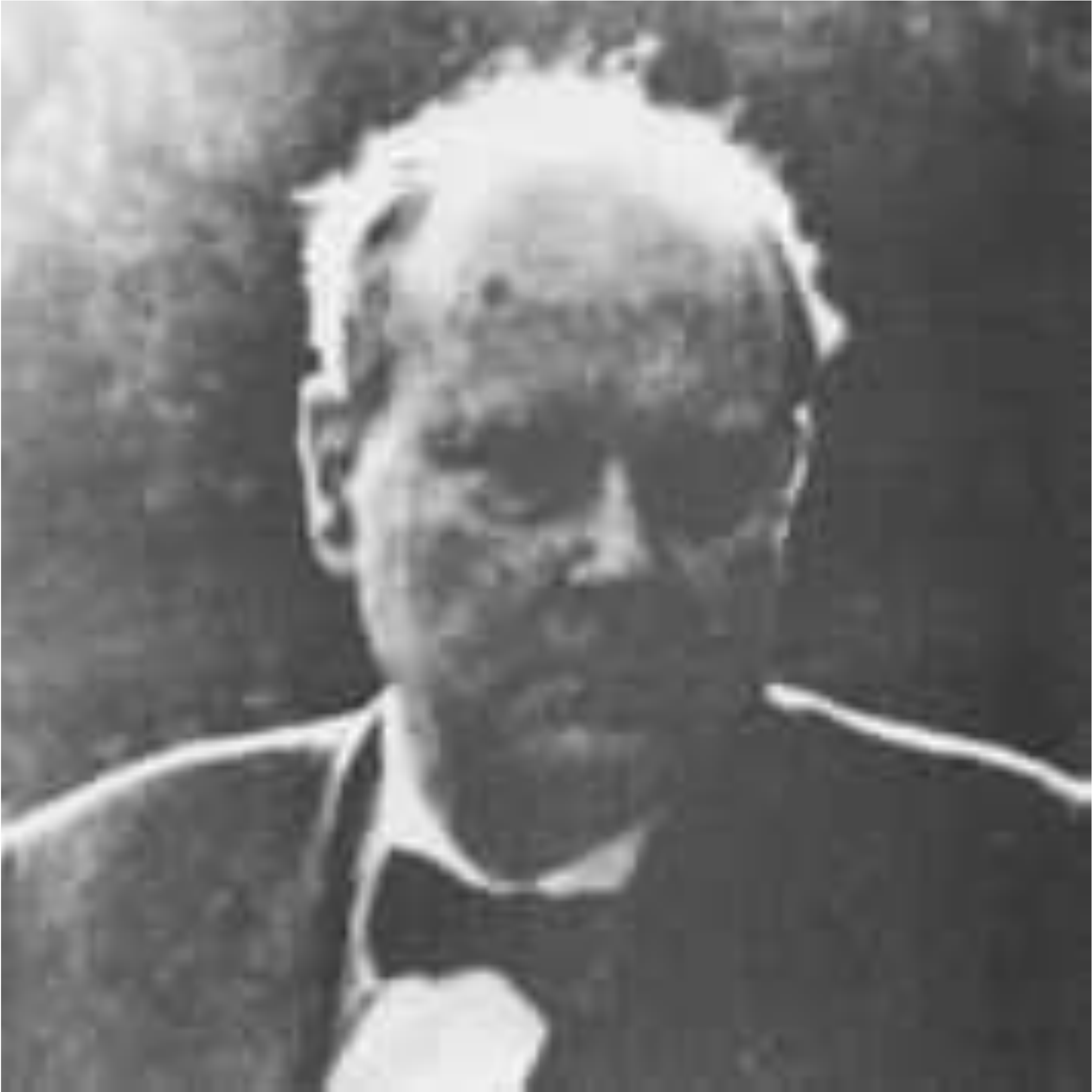}}
\newline
    {\bf Figure 1:\/} Cecil John Alvin Evelyn
\end{center}

A true ``gentleman of leisure'', he had hobbies rather than jobs. Among those hobbies was a
genuine passion for elementary geometry.  Jack and some friends,
also gentlemen of leisure, often spent time in a cafe
hand-plotting various lines and circles on large sheets of paper
in pursuit of new configuration theorems.
These plots, which today might be
routine manipulations with modern geometry software, back then
were acts of scientific inquiry.  One result of these meetings
was a self-published book ``The Seven Circles Theorem and other new theorems'',
[{\bf EMCT\/}].  (He also co-authored several papers in number theory; see
  the bibliography in [{\bf Tyr\/}].)  This book did indeed contain a
  new theorem about a configuration of touching circles, the
  Seven Circles Theorem. This theorem reads as follows:

  \begin{theorem}
    For every chain $H_1,..., H_6$ of consequently touching circles
      inscribed in and touching the unit circle
      the three main diagonals of the hexagon comprised of the points at which
      the chain touches the unit circle intersect at a common point. 
  \end{theorem}

  Figure 2 shows the Seven Circles Theorem in action.

\begin{center}
\resizebox{!}{2.5in}{\includegraphics{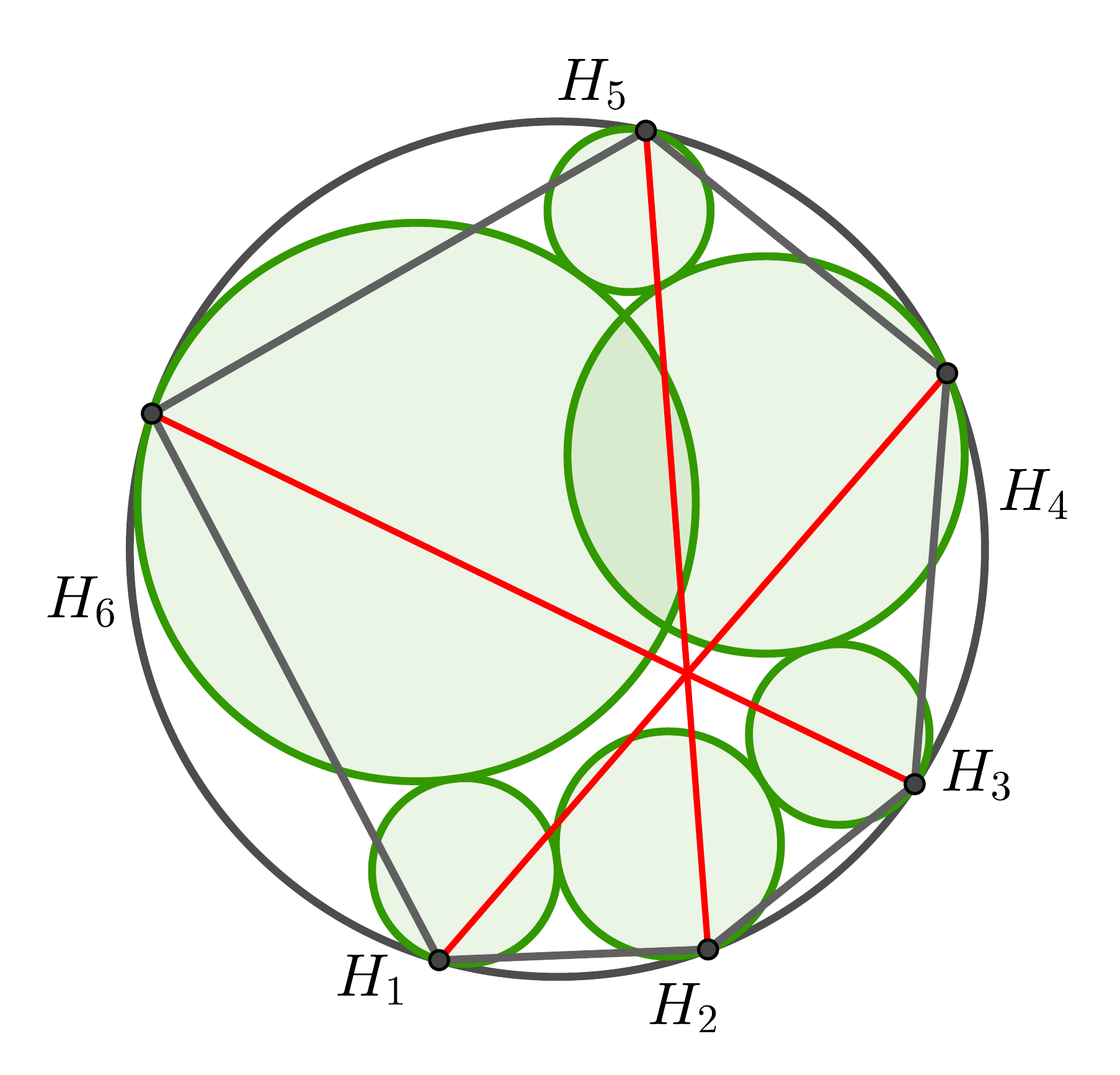}}
\newline
    {\bf Figure 2:\/} The Seven Circles Theorem
\end{center}

There are several proofs of this result.  See, for instance
[{\bf Cu\/}], [{\bf EMCT\/}], or [{\bf Ra\/}].   We noticed
that the
 Seven Circles Theorem fits naturally into the setting
 of hyperbolic geometry because everything in sight takes place inside
 the unit disk, and the open unit disk is
a common model for the hyperbolic plane.  What is interesting is
that actually the open unit disk is a model
for hyperbolic plane in {\it two\/} ways, as the Klein model and as the
Poincar\'e model.  The key to decoding the Seven Circles Theorem
is understanding the conversion between these two models.
In this note, 
we will explain the connection between the Seven Circles Theorem
and hyperbolic geometry, then prove a stronger result about
hyperbolic geometry hexagons which implies the Seven Circles
Theorem as a special case.

\section{Hyperbolic Geometry without Distances}

On the simplest level, the hyperbolic plane is a system of
points and lines which satifies the first $4$ of Euclid's
axioms and not the $5$th axiom -- the 
parallel postulate. About $180$ years ago, Bolyai,
Gauss, and Lobachevsky all discovered that one really
can make a system like this.   Here we discuss $3$ models
for the hyperbolic plane.  The first two models
are quite common and the third one is specially
introduced in order to better examine the
relationship between the first two models.
Let $\Delta$ be the open unit disk $x^2+y^2<1$  in the plane
and let $\Sigma$ denote the open northern
hemisphere in the unit sphere $x^2+y^2+z^2=1$.
\newline
\newline
    {\bf The Klein model:\/}  In this model,
the {\it points\/} are the points of $\Delta$ and    
the {\it lines\/} are the intersections
of straight lines with $\Delta$.
\newline
\newline
    {\bf The Poincar\'e model:\/} In this
    model, the {\it points\/} are the points of $\Delta$,
    and the {\it lines\/} are the intersections of circles with
    $\Delta$, provided that the circle intersects the
    boundary of $\Delta$ at right angles.
    \newline
    \newline
    {\bf The Hemisphere model:\/}
    In this model, the {\it points\/} are the points of $\Sigma$ and the
    {\it lines\/} are intersections
    of vertical planes with $\Sigma$.
    \newline

    Consider the maps
    $f: \Sigma \to \Delta$ and
    $f^{-1}: \Delta \to \Sigma$ given by the formulas
    $$f(x,y,z)=(x,y), \hskip 30 pt f^{-1}(x,y)=(x,y,\sqrt{1-x^2-y^2}).$$
    Geometrically, the map $f$ is just vertical projection.  The map $f$ carries
    lines in the Hemisphere model to lines in the Klein model, and
    the map $f^{-1}$ does the reverse.  As the
    notation suggests, $f$ and $f^{-1}$ are inverse maps. So, $f$ and $f^{-1}$
    give the conversion between the Klein and Hemisphere models.

    Now consider the maps $g: \Sigma \to \Delta$  and
    $g^{-1}: \Delta \to \Sigma$ given by the formulas
    $$g(x,y,z)=\bigg(\frac{x}{1+z},\frac{y}{1+z}\bigg),
    \hskip 15 pt
    g^{-1}(x,y)=\frac{1}{1+x^2+y^2}(2x,2y,1-x^2-y^2).$$
    Geometrically, the points $(0,0,-1)$ and $(x,y,0)$ and $g^{-1}(x,y)$ are collinear.  That is,
    we get $g(x,y,z)$ by taking the line through the south pole of the unit sphere and
    the point $(x,y,z)$ on $\Sigma$ and intersecting it with $\Delta$.  This map is known as
    {\it Stereographic Projection\/}. So, $g$ and $g^{-1}$
    give the conversion between the Poincar\'e and Hemisphere models.

    The maps $gf^{-1}$ and $fg^{-1}$ give an equivalence between the Klein and
    Poincar\'e models.  Rather than write formulas for this maps we explain what
    they do geometrically.  The map from the one model to the other just replaces
    each line of one kind with the line of the other kind that has the same ``endpoints''
    on the unit circle.  What is miraculous about this description is that there are
    underlying maps, on the level of points, which do the right things to the lines.
    So, if we have $3$ Klein-lines all containing the same {\it triple point\/}, as we do in the conclusion
    of the Seven Circles Theorem, the corresponding Poincar\'e-lines also intersect in a triple point.
    With this in mind, we show what Figure 2 looks like when we replace the three
    relevant Klein-lines with the three relevant Poincar\'e-lines.

    \begin{center}
    \resizebox{!}{2.5in}{\includegraphics{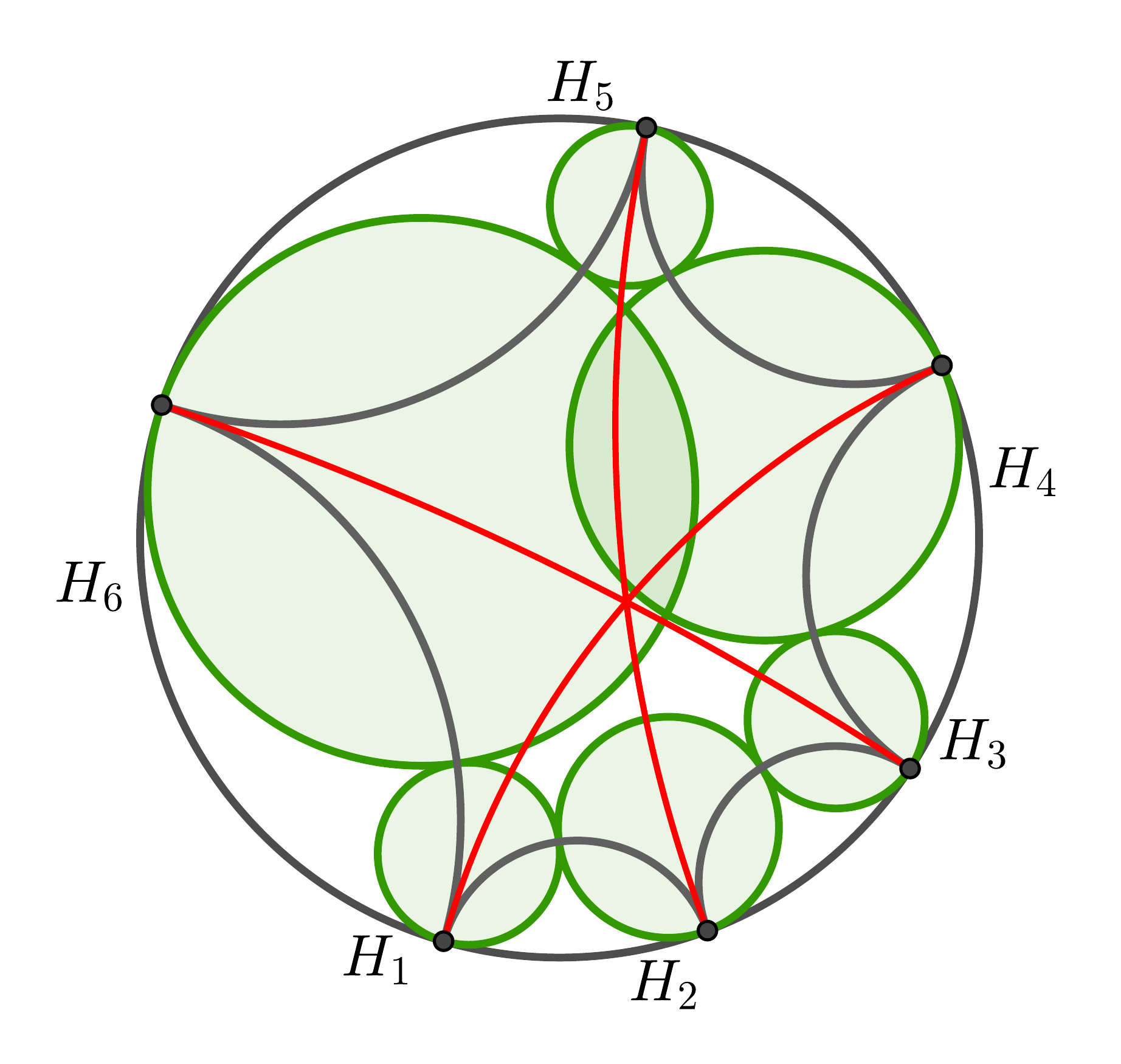}}
     \newline
    {\bf Figure 3:\/} The Seven Circles Theorem translated to the Poincar\'e model
    \end{center}

    Let $\H^2$ denote the Poincar\'e model of the hyperbolic plane.
    We hereby call Poincar\'e-lines {\it geodesics\/}, as is commonly done. Second, we denote the unit
    circle as $\partial \H^2$ and call it the {\it ideal boundary\/} of $\H^2$, as is commonly done.
    The points of $\partial \H^2$ are known as {\it ideal points\/}, even though
    technically they are not points of $\H^2$.  Finally, we say that a disk
    contained in $\H^2$, except for a single ideal point, is a {\it horodisk\/}.
    We will explain the geometric significance of such horodisks below.
    In Figure 3, the $6$ points where the horodisks intersect
    $\partial \H^2$ are ideal points.  The figure made by
    connecting the $6$ ideal points by geodesics is known as an
    {\it ideal hexagon\/}.  This ideal hexagon is shown in black.
    Here is our first (but not last) reformulation of the Seven Circles Theorem.
    
    \begin{theorem}
      \label{ss1}
      Let $P$ be an ideal hexagon.  Suppose that there are
      horodisks $H_1,...,H_6$ placed at the ideal
      vertices of $P$ in such a way that every
      two consecutive horodisks are tangent.  Then the hyperbolic
      geodesics connecting opposite vertices of $P$ meet at a triple point.
    \end{theorem}

    Once we make the conversion to the Poincar\'e model, our first reformulation is
    mostly an exercise in using new terminology.  However, after
    we go more deeply into hyperbolic geometry, we will prove a different
    reformulation that actually says something new.

    \section{Hyperbolic Geometry with Distances}
    
  Like the Euclidean plane, the hyperbolic plane is not only
    a system of points and lines but also a {\it metric space\/}.
    \newline
    \newline
    {\bf The Hyperbolic Metric:\/}
    The distance between two points $b,c \in \H^2$ is computed as follows.
    The geodesic containing $b$ and $c$ intersects $\partial \H^2$ in
    points $a,d$.  These points are labeled so that $a,b,c,d$ occur in order.
    The quantity
    \begin{equation}
      \label{dist}
      {\rm dist\/}(b,c)=\log \frac{(a-c)(b-d)}{(a-b)(c-d)}.
    \end{equation}
    is known as the {\it hyperbolic distance\/} between $b$ and $c$.
    Here we are taking advantage of the fact that we can represent
    points in the plane as complex numbers.  Thus $(a-c)(b-d)$ is
    the product of the complex numbers $a-c$ and $b-d$.  Conveniently,
    the quantity inside the $\log$ function is always a real number greater than
    or equal to $1$.
    \newline
    \newline
    {\bf Hyperbolic Isometries:\/}
    Just as the Euclidean distance function, typically defined in terms of
    sums of squares, exhibits a surprising rotational symmetry, so does the
    hyperbolic distance function.  A {\it linear fractional transformation\/} of
    the complex plane $\C$ is a map of the form
    \begin{equation}
      T(z)=\frac{az+b}{cz+d}, \hskip 30 pt
      ad-bc \not =0
    \end{equation}
    Such a map is called a {\it hyperbolic isometry\/} when both $T$ and $T^{-1}$
    map $\H^2$ to itself.  The reason for the name, in this case, is that
    $T$ preserves hyperbolic distances:  ${\rm dist\/}(T(b),T(c))={\rm dist\/}(b,c)$ for all points
    $b,c \in \H^2$.  One can verify this claim with a straightforward
    algebraic manipulation.   One can map any point of $\H^2$ to any other
    point by such a hyperbolic isometry.  For instance, the map
    $$T(z)=\frac{z-r}{rz-1}$$
    is a hyperbolic isometry which has the property that $T(\pm 1)=\pm 1$ and $T(r)=0$.
    The rotations about
    the origin are also hyperbolic isometries.  Using these two kinds of
    maps, you can probably convince yourself of our claim that one can map
    any point of $\H^2$ to any other point using a hyperbolic isometry.
    \newline
    \newline
    {\bf Hyperbolic Disks:\/}
    A {\it hyperbolic disk\/} in $\H^2$ any set of the form
    $$B(z_0,r)=\{z|\ {\rm dist\/}(z,z_0) \leq r\}.$$
    Here $z_0$ is a point of $\H^2$, called the {\it hyperbolic center\/} of the hyperbolic disk.
    The boundary of the hyperbolic disk is called a {\it hyperbolic circle\/}.
    The boundary of a horodisk is a {\it horocircle\/}.
    
    The hyperbolic disks
    centered at the origin are round Euclidean disks, by symmetry.  Also, the
    hyperbolic isometries map Euclidean disks to Euclidean disks.  Since these
    hyperbolic isometries also preserve the hyperbolic distance, we see that
    as sets the hyperbolic disks and circles are exactly the same as the
    Euclidean disks and circles contained in $\H^2$.  Note, however, that
    the hyperbolic
    and Euclidean center of a disk are typically distinct, as are the
    hyperbolic and Euclidean radii.  If we keep the Euclidean size of a disk
    the same and start shifting it over until it becomes tangent to $\partial \H^2$,
    the hyperbolic center moves to the tangency point and the radius tends to $\infty$.
    So, one can view a horodisk as a disk of infinite radius centered at an ideal point.
    The first two pictures in Figure 4 show hyperbolically concentric circles, and
    the last picture shows hyperbolically concentric horocircles.
        
    \begin{center}
    \resizebox{!}{1.3in}{\includegraphics{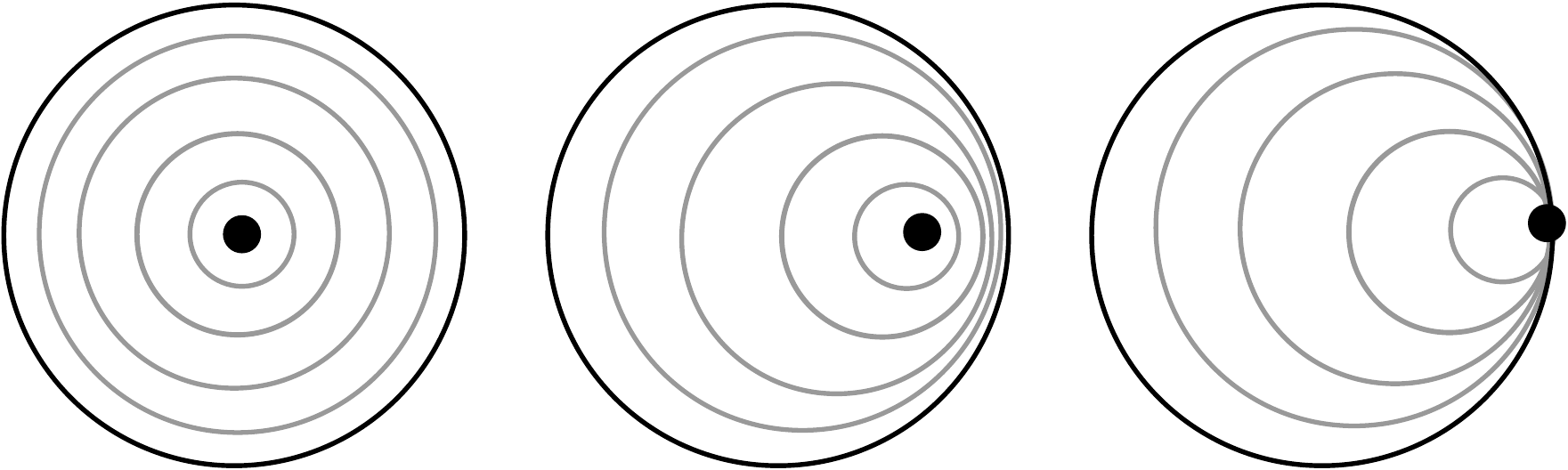}}
     \newline
         {\bf Figure 4:\/} Hyperbolically concentric circles and horocircles.
    \end{center}

    \noindent
        {\bf Alternating Perimeter:\/}
        A {\it hyperbolic polygon\/} has the same definition as a Euclidean polygon
        except that everything takes place in $\H^2$.  Namely, it is a closed loop
        made by connecting together finitely many hyperbolic geodesic segments.
        The endpoints of these segments are the {\it vertices\/} and the segments
        themselves are the {\it edges\/}.   A hyperbolic polygon always has a
        perimeter, namely the sum of the hyperbolic distances between the vertices.
        Equivalently, we could say that the perimeter is the sum of the lengths of the edges.
        When the polygon has an even number of sides, it also has an {\it alternating perimeter\/}.
        This is defined to be alternating sum of the lengths of the edges, namely
        $(S_1+S_3+S_5+...)-(S_2+S_4+S_6+...)$.

        An {\it ideal polygon\/} has the same definition as above, except that the
        vertices are all ideal points.  Surprisingly, the alternating
        perimeter of an ideal polygon makes sense, even though the individual
        terms in the sum are infinite.  To see this for hexagons (which is the case we
        care about) imagine that we have a sequence $\{P_n\}$ of hyperbolic hexagons
        converging to some ideal hexagon $P_{\infty}$.
        We place disks at the vertices of $P_n$, in such a way
        that consecutive ones do not overlap, and we
        compute the alternating perimeter instead as
        \begin{equation}
          \label{horosum}
          (S_1'+S_3'+S_5')-(S_2'+S_4'+S_6')
          \end{equation}
        where $S_k'$ is the length of the portion of the $k$th side of $P_n$ that
        lies outside the two disks centered at its endpoints. If we replace one of
        our disks by another one with the same center, we are adding some amount
        to some $S'_k$ and subtracting the same amount to $S'_{k+1}$.  So, the
        modified sum does not depend on which disks we choose.   Taking a
        limit as $n \to \infty$, we see that the alternating perimeter of
        $P_{\infty}$ can be defined as the same kind of sum as in
        Equation \ref{horosum}, except that $S_k'$ denotes the (finite!) length
        of the portion of the $k$th side of $P_{\infty}$ that lies outside
        the two horodisks centered at its endpoints.

        \subsection{The Main Result}

        An ideal hexagon $P$ determines a
        small triangle $T_P$, the geodesic triangle bounded by the
        three geodesics connecting opposite sides of $P$.  In Figure 5 below,
        the red triangle is $T_P$.  Here is the main result.

        \begin{theorem}
         \label{ss3}
         For any ideal hexagon $P$, the alternating perimeter of $P$ is, up to sign, twice the
         perimeter of $T_P$.
        \end{theorem}
        
        \startproof
        Define a {\it semi-ideal triangle\/} to be a hyperbolic triangle with
        two ideal vertices and one vertex in $\H^2$.  The yellow triangles $Y_1,Y_2,Y_3$
        on the left side in Figure 5 are semi-ideal triangles.  The green
        triangles $G_1,G_2,G_3$ on the right side in Figure 5 are also semi-ideal triangles.
        Note that the green triangles overlap. They all contain $T_P$.

        \begin{center}
        \resizebox{!}{1.9in}{\includegraphics{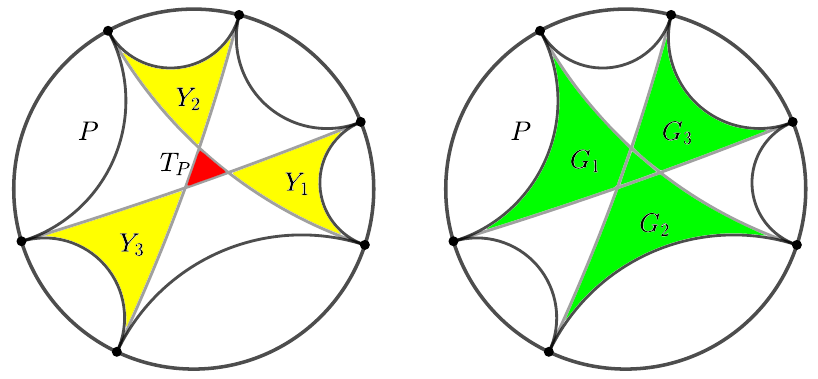}}
         \newline
         {\bf Figure 5:\/} The ideal hexagon $P$ and the small triangle $T_P$.
        \end{center}

        Let $L_1,L_2,L_3$ denote the sides of a semi-ideal triangle $V$, with the convention
        that $L_3$ connects the two ideal vertices of $V$.  We delete disjoint horodisks from the
        two ideal vertices of $V$ and define
        \begin{equation}
          A(V)=L'_1+L'_2-L'_3,
        \end{equation}
        Where $L'_j$ is the length of the portion of $L_j$ outside the two horodisks.
        This definition is very similar to the definition of the alternating perimeter
        above, and it does not depend on which horodisks we remove.

        Note that $G_k$ and $Y_k$ share a vertex for $k=1,2,3$.
        We have a hyperbolic isometry $I_k$ such
        that $I_k(Y_k)=G_k$.  This is most easily seen if we normalize the picture so that
        the vertex $Y_k \cap G_k$ is the Euclidean origin; in this case $I_K$ is just reflection in the origin.
        Hence $A(Y_k)=A(G_k)$.  Summing up, we have
        \begin{equation}
          \label{trisum}
          A(Y_1)+A(Y_2)+A(Y_3)-
          (A(G_1)+A(G_2)+A(G_3))=0.
          \end{equation}
        The terms in the sum on the left side
        of Equation \ref{trisum} can be divided into two types:
        Those which come from geodesics connecting ideal points,
        and the rest.  When we sum up the terms of the first kind,
        we get $A(P)$ (up to sign) When we sum up the rest of the terms,
        we get twice the perimeter of $T_P$ because, so to speak,
        the lengths of the rest of the geodesics bounding $T_P$
        are counted in pairs with opposite signs.
        \endproof

        Now that we have Theorem \ref{ss3}, we can give our final
        reformulation of the Seven Circles Theorem.
        Say that a hexagon $P$ has {\it point reflection symmetry\/} if
        there is a nontrivial hyperbolic isometry $T$, which has a single
        fixed point, such that $T(P)=P$ and $T^2$ is
        the identity.  We call $T$ a {\it point reflection\/}.
        The map $T$ swaps the opposite vertices of $P$.
        Here is our final reformulation of the Seven Circles Theorem.

        \begin{theorem}
          \label{ss2}
          The following are equivalent for an ideal hexagon $P$.
          \begin{enumerate}
          \item The geodesics connecting opposite vertices of $P$ meet at a triple point.
          \item $P$ has point reflection symmetry.
          \item The alternating perimeter of $P$ is $0$.
          \end{enumerate}
         \end{theorem}

        \startproof
        ($1 \to 2$): Suppose that the first condition holds.
        The conditions above are unchanged if we move the picture by a
        hyperbolic isometry.  We do this in such a way that the
        triple point lies at the origin of the unit disk.   But then, in the Poincar\'e model,
        the three geodesics connecting
        the opposite vertices of $P$
        are all Euclidean diameters of the unit disk.  This means that the Euclidean
        symmetry $z \to -z$ is a symmetry of $P$.
        But this particular symmetry is also a hyperbolic
        isometry.  Hence $P$ has point reflection symmetry.

        ($2 \to 3$):  We choose the horidisks
        $H_1,H_2,H_3$ arbitrarily, then afterwards we set $H_{k+3}=T(H_k)$ for $k=1,2,3$.   This makes all of
        Figure 3 invariant under the point reflection $T$.  But then, referring to Equation
        \ref{horosum}, we have $S'_{k+3}=S'_k$ for $k=1,2,3$.  Hence the sum in
        Equation \ref{horosum} is $0$.

        ($3 \to 1$):  If $P$ has alternating perimeter $0$ then, by Theorem
        \ref{ss3}, the triangle $T_P$ has perimeter $0$.  This means that
        $T_P$ is actually a single point.
        \endproof

        If we have an ideal hexagon as in
        Theorem \ref{ss1}, then its alternating perimeter is $0$ because in Equation \ref{horosum}
        we gave $S_1'=...=S_6'=0$.  So, Theorem \ref{ss2} implies Theorem \ref{ss1}, and hence the
        Seven Circles theorem.  Theorem \ref{ss2} also reveals that
        all the instances of the Seven Circles Theorem involve ideal hexagons having point reflection
        symmetry, and this makes the Seven Circles Theorem obvious.

\hskip 5pt

\noindent
\textbf{Acknowledgments.} We are grateful to the Mathematics Institute of the Heidelberg University for providing a stimulating environment, and we are especially grateful to Sergei Tabachnikov for bringing the Seven Circles Theorem to our attention and for inspiring discussions on the topic.

\section*{References}

\noindent
[{\bf Cu\/}] H.\, M.~Cundy, {\it The Seven-circles Theorem\/}, Math. Gazette. {\bf 62\/} (1978) pp 200-203.
\vskip 7 pt

\noindent
[{\bf EMCT\/}] C.\,J.\,A.~Evelyn, G.\,B.~Money-Coutts, J.\,A.~Tyrrell, {\it The seven circles theorem and other new theorems\/}, Stacey International, London, 1974.
\vskip 7 pt

\noindent
[{\bf Tyr\/}] J.\,A.~Tyrrell, {\it Cecil John Alvin Evelyn\/}, Bull. London Math. Soc. {\bf 9\/} (1977) pp 328-329.
\vskip 7 pt

\noindent
[{\bf Ra\/}] S.\, Rabinowitz, {\it The Seven Circles Theorem\/}, Pi Mu Epsilon Journal {\bf 8\/} (1987) pp 441-449.
\end{document}